\documentclass[a4paper,oneside,11pt]{article}

\usepackage[francais]{babel} 
\usepackage{vmargin}
\usepackage{epsfig}
\psfigdriver{dvips}
\usepackage{graphicx}
\usepackage{color}
\usepackage[latin1]{inputenc}
\usepackage{amsmath}
\usepackage{amssymb}
\usepackage[all]{xy}

\newcommand{\M}{{\mathbf M}}

\newcommand{\Ext}{{\rm Ext}}
\newcommand{\Q}{{\mathbf Q}}
\newcommand{\im}{{\rm Im}}

\newcommand{\MF}{{\mathcal MF}}

 \newcommand{\adj}{{\rm adj}}
 \newcommand{\cok}{{\rm cok}}
 \newcommand{\coker}{{\rm coker}}
 \newcommand{\GL}{{\rm GL}}
  \newcommand{\Hom}{{\rm Hom}}
  \newcommand{\Def}{{\rm Def}}

\begin{document}

\newtheorem{thm}{Théorème}
\newtheorem{cor}[thm]{Corollaire}
\newtheorem{lem}[thm]{Lemme}
\newtheorem{prop}[thm]{Proposition}
\newtheorem{dfn}[thm]{Définition}
\newtheorem{rem}[thm]{Remarque}
\newtheorem{exm}[thm]{Exemple}
\newcommand{\al}{\alpha}
\newcommand{\A}{\mathcal{A}}
\newcommand{\Aa}{\mathcal{A}^\alpha}
\newcommand{\ta}{\theta}

\sloppy 

\begin{center}
{\Large\bf  P\'eriodicit\'e de Kn\"orrer \'etendue}
\end{center} 
\vspace{0,5cm}
\begin{center}
 Jos\'e Bertin, \quad Fabrice Rosay
\end{center}
\begin{center}
\small { Institut Fourier,  
 jose.bertin@ujf-grenoble.fr, frosay@ujf-grenoble.fr}
\end{center}

\begin{abstract}
 Nous  \'etudions la th\'eorie des d\'eformations des Factorisations Matricielles, eventuellement munies d'une structure orthogonale ou symplectique. Nous discutons  et g\'en\'eralisons  dans diff\'erents contextes les th\'eor\`emes de p\'eriodicit\'e de Kn\"orrer et Hori-Walcher.\\
-----------------------------------
 We study the deformation theory aspects of Matricial Factorizations, possibly with an orthogonal or symplectic structure. We discuss and extend the Kn\"orrer  and Hori-Walcher periodicity theorems. 
  \end{abstract}
  \bigskip
  \bigskip
  \emph{\underline{\textbf  English version (abridged)}}
  \bigskip
  \bigskip
  
  In this note  our aim is to introduce orthogonal and symplectic matrix factorizations of a polynomial $w$, and their twisted  counterpart, motivated by Faltings paper \cite{fal}. Our main result is a Kn\"orrer type periodicity theorem in the spirit of \cite{orii}.  Let us fix an algebraically closed field  $k$ of   characteristic zero. Throughout a ring is a complete local $k$-algebra, with  residue field $k$. Let  $w(x)\in k[x_1,\cdots,x_n], \, w(0) = 0$ be a polynomial defining an hypersurface with isolated singularities (the potential). For deformation-theoretic reasons we will see $w$    in $R = k[[x_1,\cdots,x_n]]$. A matrix factorization (MF in short) of $w$ over a ring $A$ \cite{ori}, \cite{Kn}, is the data of a 
  $\mathbb Z/2\mathbb Z$-gradued   free $R\hat\otimes A$-module $\M  = M_0\oplus M_1$, together with an odd map $\Q: \M \to \M$, such that $\Q^2 = w.1$. If $A = k$, the subscript $k$ will be dropped.  Define $\M[1]$ as the  shift graded module with odd operator $-\Q[1]$ \  \cite{ori},\cite{KoRo}. Let be given $\M, \, \M'$  two  MF over $A$, then the  graded module $\hom_{R\hat\otimes A}(\M,\M')$ is a graded dg-module with differential $D(f) = [\Q,f] = \Q'f-(-1)^{\vert f\vert} f\Q$ the graded commutator. The degree zero MF-morphisms,  in short the morphisms, are those $f\in \hom^0(\M,\M')$ with $D(f) = 0$. Define    ${\rm Ext}(\M,\M') = \ker Q/\im Q$ \cite{KoRo}. It is a graded finitely generated $A$-module. The last  construction we need is the tensor product $\M\otimes \M'$. It is simply the graded tensor product  of the corresponding dg-modules.  Finally we must observe  the category of MF (over $A$) of $w$ is a DG-category ${\MF}_A(w)$, which yields  if we take morphisms up homotopy the triangulated category ${\rm MF}_A(w)$ \cite{KoRo}, \cite{Or}.
  
  Let $\M$ be a MF over $A$. Define $\M^T$ (resp. $\M^*$) as the  graded dual module together with the transpose of $\Q$ (resp. $\M^T[1]$). A quadratic (resp. symplectic) form on $\M$ is an isomorphism $b: \M \stackrel\sim \to \M^*$ such that $^tb = -b$ (resp. $^tb = b$), where  $^tb$ denote the transpose of $b$.
  
  A twisted quadratic (resp. symplectic) form on $\M$ is an isomorphism $q: \M \stackrel\sim \to \M^T$ such that $^tq = q$ (resp. $^tq = -q$).  Let us denote $\MF^{\pm} $ (resp. $\MF^{\omega,\pm}$) the category of MF equipped with a quadratic (+), or symplectic (-) (resp. twisted quadratic or symplectic) form, the morphisms being  the morphisms of underlying MF.
  
  If $\M$ and $\M'$ are equipped with a quadratic (resp. symplectic, twisted quadratic,.....) $b$ (resp. $q$),   the adjoint morphism  of $f: \M \to \M'$ is defined, yielding an  operator $\adj: \hom(\M,\M') \to \hom(\M',\M)$, as well $\adj: {\rm Ext}(\M,\M') \to {\rm Ext}(\M',\M)$. 
  The tensor product  previously defined extends  to (twisted) quadratic or symplectic MF, following the rule $b\otimes q'$  or $q\otimes b'$  is a quadratic (or symplectic) form on $\M\otimes \M'$, and $b\otimes b'$, or $q\otimes q'$ is a twisted quadratic (or symplectic) form  (Prop. 3).
  
The matrix factorizations have good deformation properties. Let   $\overline\M$ be a MF over $k$. There is a well-defined deformation functor $\Def_{\overline\M}$ (resp. a restricted functor keeping $w\in R$ fixed) with tangent space $t_{\overline\M}$ which fits into an exact sequence  
\begin{displaymath} 0\to \Ext^1(\overline\M,\overline\M) \to t_{\overline\M} \to \mathcal J\to 0\end{displaymath}
(resp. $t_{\overline\M} = \Ext^1(\overline\M,\overline\M)$), where $\mathcal J$ denotes the subspace of $Q$-exact functions in the Tjurina ring. The obstructions to an infinitesimal lifting lie in the cokernel of the natural map $\mathcal J \to \Ext^0(\overline\M,\overline\M)$ (resp. in $\Ext^0(\overline\M,\overline\M)$).  Similarly if $\overline\M$ is a (twisted) quadratic or symplectic MF, there is a corresponding deformation functor, with analogous description with however $\Ext^k(\overline\M,\overline\M)^\pm$  the subgroup of selfadjoint, or anti-selfadjoint elements in place of  $\Ext^k(\overline\M,\overline\M)$ (Prop. 4).
  
  Let us denote $(x,y)$ the MF $k[[x,y]]\stackrel x\to k[[x,y]]\stackrel y \to k[[x,y]]$ of $w = xy$. Let $(P,Q)$ be a MF  with entries in $A$, and potential $\pi\in {\cal M}_A$.  The Faltings-Kn\"orrer functor $\theta: \MF_A(\pi)\to \MF_{A[[x,y]]}(xy-\pi)$ is given by the tensor product $\theta (P,Q) = (Q,-P)\otimes (x,y)$.
If we endow $(x,y)$ with its essentially unique quadratic form, the functor $\theta$ extends   to 
\begin{displaymath}\theta: \MF^{\omega,\pm} (\pi) \longrightarrow \MF^{\pm}(xy-\pi)\end{displaymath}
We then focus on the deformation-theoretic properties of $\theta$.   Our result is the fact the Faltings-Kn\"orrer functor yields an equivalence of categories both at $\MF$ and MF levels  (Thm 5). If we restricts to quadratic  or symplectic MF, the square $\theta^2$ yields an equivalence  between appropriated categories (Th\'eor\`eme 6), result  in the spirit of \cite{orii}.

\section{Introduction}

 Dans cette note   notre objectif est de traduire dans le cadre des factorisations matricielles  quelques r\'esultats de Faltings \cite{fal} portant  sur la structure locale   des espaces de modules de $G$-fibr\'es vectoriels sur une courbe alg\'ebrique semi-stable, $G$ \'etant l'un des groupes,    lin\'eaire,   orthogonal ou   symplectique.  On sait que les factorisations matricielles du polyn\^ome $w(x)\in k[x_1,\cdots,x_n]$ correspondent de mani\`ere essentiellement bijective aux modules de Cohen-Macaulay maximaux sur l'hypersurface affine  $w = 0$, suppos\'ee  avec des singularit\'es  isol\'ees \cite{Kn}.  Ceux de ces modules $E$ qui sont   bilin\'eaires,  quadratiques ou symplectiques, correspondent \`a des factorisations matricielles appel\'ees dans le pr\'esent article, quadratiques ou symplectiques. Une version tordue est obtenue en demandant que le module $E$ soit isomorphe au dual du module de syzygy $\Omega(E)$. Le foncteur de Faltings-Kn\"orrer \'echange ces deux situations \cite{fal}, \cite{Kn}, \cite{Or}.   On montrera que   le carr\'e du foncteur de Kn\"orrer  induit   une \'equivalence  entre  factorisations  matricielles bilin\'eaires  d\'efinies sur un anneau local noetherien complet $R$, et celles sur $R[[x,y,u,v]]$. Dans cette \'equivalence,  factorisations matricielles orthogonales et   symplectiques sont \'echang\'ees. La situation est analogue \`a celle \'etudi\'ee   dans un contexte diff\'erent\footnote{Ce contexte ''Orientifold'' suppose que le ''potentiel'' $w(x)$ est anti-invariant pour une involution $\tau$ de $k[x_1,\cdots,x_n]$.}par  Hori et Walcher \cite{orii}. Les r\'esultats d\'ecoulent    de la description des d\'eformations  verselles des factorisations matricielles  de $xy  $.

\section{ Cat\'egories de Factorisations matricielles}

Dans ce texte $k$ est un corps  alg\'ebriquement clos de caract\'eristique nulle. Un anneau  signifie une $k$-alg\`ebre locale noetherienne compl\`ete, de corps r\'esiduel $k$. On fixe un polyn\^ome $w(x)\in k[x_1,\cdots,x_n]$, sans terme constant, d\'efinissant une hypersurface avec  singularit\'ees isol\'ees.  On regardera   $w(x) $ comme \'el\'ement de $R = k[[x_1,\cdots,x_n]]$.  
  
  Une factorisation matricielle  \cite{orii},  \cite{Kn}  (MF en abr\'eg\'e) de rang $r\geq 1$ de $w$ \`a coefficients dans l'anneau $A$, est la donn\'ee d'un couple de matrices  $\M = (\varphi, \psi)\in M_r(R\hat\otimes A)$  avec $\varphi\psi = \psi\varphi = w.1_r$, et  $\varphi (0) = \psi (0) = 0$.   
(Voir \cite{KoRo} pour une d\'efinition plus g\'en\'erale). Dans la suite $w$ est fix\'e\footnote{Dans la section 4, on  d\'eformera  $w$, qui pourra \^etre plus g\'en\'eralement un \'el\'ement de l'id\'eal maximal de $R\hat\otimes A$ relevant $w$.}. On peut voir   $\M  $ comme un diagramme de $R_A$-modules libres de rang $r$ (le rang de la MF)
 $ M_0 \stackrel\varphi\longrightarrow M_1 \stackrel\psi\longrightarrow M_0 \quad(M_0 = M_1 = R_A^r). $
   
Un morphisme $f = (S,T): \M = (\varphi,\psi) \rightarrow \M' = (\varphi',\psi')$ est un morphisme   de diagrammes, soit un couple de matrices $(S,T\in M_r(R\hat\otimes A))$ tel que 
$\varphi'  S = T\varphi,\, \psi'T =S\psi$. Si $S$ et $T$ sont inversibles  $f$ est un isomorphisme.  Le groupe $\GL_r\times \GL_r$ agit sur les MF  par $(S,T).(\varphi,\psi) = (T\varphi S^{-1},S\psi T^{-1})$.  Il est utile de voir $\M$ comme un objet $\mathbb Z/2\mathbb Z$-gradu\'e $\M  = M_0\oplus M_1$.     On d\'efinit la MF translat\'ee $\M[1]$ comme  \'etant 
 $M_1 \stackrel{- \psi}\longrightarrow M_0 \stackrel{-\varphi}\longrightarrow M_1.  $
Soit l'op\'erateur impair $\Q: \M \to \M$:
\begin{equation}\label{operator}
\mathbf{Q} = \left(\begin{array}{cc}0&\psi\\ \varphi&0\end{array}\right): \M \rightarrow \M \quad (\mathbf{Q}^2 = w.1)  \end{equation}
Si $\M, \, \M'$ sont deux MF  de $w$ \`a coefficients dans $A$, le groupe   $\hom(\M,\M') = \hom^0(\M,\M')\oplus \hom^1(\M,\M')$ est  un $\mathbb Z/2\mathbb Z$- groupe diff\'erentiel gradu\'e, de diff\'erentielle
 $D(f) = Q'f - (-1)^{\vert f\vert} fQ $  (on note $\vert f \vert$ le degr\'e de $f$).
On pose 
\begin{equation}\label{ext}\Ext^\bullet(\M,\M') = \ker Q/ \im Q = \Ext^0(\M,\M')\oplus \Ext^1(\M,\M')\end{equation}
Comme $w(x)$ n'a que des singularit\'es isol\'ees $\Ext^\bullet(\M,\M')$ est un $A$-module de type fini.   Le produit tensoriel $\M\otimes \M'$ est la MF  qui a pour op\'erateurs
\begin{equation}\label{tens}\Phi =  \left(\begin{array}{cc}1\otimes \varphi'&\psi\otimes 1\\ \varphi\otimes 1&-1\otimes \psi'\end{array}\right),\quad \Psi =  \left(\begin{array}{cc}1\otimes \psi'&\psi\otimes 1\\ \varphi\otimes 1&-1\otimes \varphi'\end{array}\right)\end{equation}
 En r\'esum\'e,  
les  factorisations matricielles de $w(x)$ \`a coefficients dans $A$ forment une DG-cat\'egorie    $\MF_A$ avec c\^ones,   et  translation  involutive $\M \to \M[1]$ \cite{Or}. La cat\'egorie triangul\'ee  associ\'ee ${\rm MF}_A$ a  pour objets   les MF, et pour morphismes $\M \to \M'$ (de degr\'e z\'ero),  les \'el\'ements de $\Ext^0(\M,\M')$. 
 
\emph{\underline{ Exemple 1}:} Si $R = k[[x,y]]$, et $w = xy$, toute MF de $w$ \`a coefficients dans $A = k$, de rang $d$, est isomorphe \`a  $\M_{p,q}  \,\, (p,q\geq 0, \, p+q = d)$, avec pour op\'erateurs
  $\varphi = \left(\begin{array}{cc} x1_p & 0\\0 & y1_q\end{array}\right), \quad \left(\begin{array}{cc}y1_p& 0\\  0&x1_q\end{array}\right)$, 
 soit $\M_{p,q} =   \M_{1,0}^{\otimes p}\otimes   {\M}_{1,0}[1]^{\otimes q} $.

\section{ Factorisations orthogonales et symplectiques}

 Soit $\M: \,  M_0 \stackrel\varphi\longrightarrow M_1 \stackrel\psi\longrightarrow M_0 $ une MF de rang $d$.  Posons    $\M^T: \,  M_0^* \stackrel {^t\psi}\longrightarrow M_1^* \stackrel {^t\varphi}\longrightarrow M_0^*$, et $\M^*: \,  M_1^* \stackrel{-{^t\varphi}}\longrightarrow M_0^* \stackrel {-{^t\psi}}\longrightarrow M_1^*$. On a $\M^* = (\M^T)[1] = (\M[1])^T$.  Ces op\'erations s'\'etendent en des foncteurs anti-involutifs $\MF_A \longrightarrow \MF_A$.
  On a des identifications canoniques
 $ (\M\otimes \M')^T \cong \M^T\otimes \M'^T\cong \M^*\otimes \M'^*$, et $  
 (\M\otimes \M')^* \cong { \M}^T\otimes \M'^* \cong \M^* \otimes {\M'}^T$ 
    
 \begin{dfn} On appelle forme quadratique (resp. symplectique) sur $\M$, la donn\'ee d'un isomorphisme  dans $\MF_A$, 
 $b: \M\stackrel\sim \to \M^*$ (morphisme de degr\'e un   $\M\stackrel\sim \to \M^T$) tel que ${^tb} = -b$ (resp. $b$), et $Q^{\adj} := b^{-1}{^tQ}b = -Q$.
 On appelle forme quadratique (resp. symplectique) tordue sur $\M$, la donn\'ee d'un isomorphisme $q: \M \stackrel\sim\to \M^T$ tel que $Q^{\adj} = Q$, et ${^tq} = q$ (resp. $-q$).\end{dfn}
 D'une autre mani\`ere si $b = \left(\begin{array}{cc} 0 & b_1\\b_0 & 0\end{array}\right), \, \, b_0: M_0\to M_1^*, \, b_1: M_1\to M_0^*$, alors ${^tb}_0 =  -b_1$ (resp. $b_1$), et $\varphi^{\adj} := b_1^{-1}{^t\varphi}b_0 = -\varphi, \, \psi^{\rm adj} = -\psi$.  Explicitation analogue dans le cas tordu. On note $\MF_A^{\pm}$ (resp. $\MF_A^{\omega,\pm}$) la cat\'egorie dont les objets sont les MF avec structure orthogonale (+), ou  symplectique ($-$), (resp. orthogonale tordue, symplectique tordue).     Si $A = k$, l'indice $A$ sera omis. Un isomorphisme $f: (\M,b)\stackrel \sim\to (\M',b')$ est un isomorphisme de MF  tel que $f = b^{-1}{^tf}b'$, d\'efinition analogue  dans le cas tordu. Le groupe de jauge agit sur les  MF (orthogonales, ....) vues comme des couples $(Q,b)$ selon la r\`egle $F(Q,b) = (F^{-1}QF , {^tF}bF)$.

   Posons $ \cok (\M) := \coker (\varphi)$. Le foncteur $\cok$ \'etablit une \'equivalence de cat\'egories de ${\rm MF}_k$ sur la cat\'egorie (stable) des $R/(w)$-modules de Cohen-Macaulay maximaux \cite{Kn}. Le module des syzygies $\Omega(\cok (\M))$ est $\cok (\M[1]) = \coker (\psi)$. On a   $\cok (\M^*) = \cok(\M)^*$, et $\cok (\M^T) = \Omega(\cok(\M))^*$. Une forme quadratique (resp. symplectique) sur $\M$ induit une forme quadratique (symplectique) sur $ \cok (\M)$. Une forme quadratique (resp. symplectique) tordue induit  un isomorphisme $c: \cok (\M) \stackrel \sim \to (\Omega(\cok (\M))^*$ tel que ${^t\Omega}(c) = \pm c$.  R\'eciproquement toute structure orthogonale (resp. ....)  sur un module de Cohen-Macaulay maximal provient d'une MF orthogonale (resp.  ...).

 Si $\M, \, \M'$ sont deux objets de $\MF^+$ (resp. $\MF^-$,....),  l'application d'adjonction (de degr\'e z\'ero) est
  $\adj:  \Hom^\bullet_{\MF^\pm}(\M,\M') \longrightarrow \Hom^\bullet_{\MF^\pm}(\M',\M), \quad f\mapsto f^{\adj}=b^{-1}{^tf}b $
(resp. $\adj:  \Hom^\bullet_{\MF^{\omega,\pm}}(\M,\M') \longrightarrow \Hom^\bullet_{\MF^{\omega,\pm}}(\M',\M)$). 

 \begin{lem}  On a la r\`egle de commutation
  \begin{equation}\label{com} D(f)^{\adj} = (\pm)(-1)^{\vert f\vert} D(f^{\adj}) \end{equation}
 avec $+$ dans le cas orthogonal ou symplectique, et $-$  dans le cas orthogonal ou symplectique tordu.\end{lem}
  L'application d'adjonction est d\'efinie au niveau de la cat\'egorie  homotopique ${\rm MF}_A$, et conduit \`a   $\adj: \Ext^\bullet(\M,\M') \to \Ext^\bullet(\M',\M)$.  
   On utilisera le signe $\epsilon = +1$ (resp. $\epsilon = -1$) pour d\'esigner un type de structure, une structure orthogonale ou orthogonale tordue (resp. symplectique, symplectique tordue)
 \begin{prop} Soit $(b,\epsilon)$   ou $(q,\epsilon)$ (resp. $ (b',\epsilon'),(q',\epsilon')$)   une structure bilin\'eaire ou  bilin\'eaire tordue sur $\M$ (resp.  sur $\M'$),  alors $b\otimes q'$ ainsi  que  $q\otimes b'$ 
 d\'efinit une structure bilin\'eaire de type $\epsilon\epsilon'$. Par contre $q\otimes q'$ (resp. $b\otimes b'$)  d\'efinit une structure bilin\'eaire tordue de type $\epsilon \epsilon'$ (resp. $-\epsilon \epsilon'$).\end{prop}
 
 \emph{\underline{ Exemple 2}:} Soient des entiers $h_i\geq 2 \, (1\leq i\leq d)$. Pour $n_i\in [1,h_i-1]$, soit $\M_{n_i}$ la MF de rang un d\'efinie sur $k[[x_i]]$ par $(x_i^{n_i},x_i^{h_i-n_i})$. Soit au-dessus de $k[[x_1,\cdots,x_d]]$ la MF $\M_{n_1,\cdots,n_d}= \bigotimes_{i=1}^d \M_{n_i}$. Il y a sur  $\M_{n_1,\cdots,n_d}$ une (essentiellement) unique structure bilin\'eaire  si $d$ est impair (resp. bilin\'eaire tordue si $d$ est pair), qui est orthogonale ou symplectique selon la parit\'e de $m = [d/2]$.

 \section{D\'eformations et foncteur de Faltings-Kn\"orrer}

 Fixons   $\overline{\M} = (\overline{\varphi},\overline\psi)$ une MF de rang $d$, d\'efinie  sur $k$. On suppose que le potentiel $w(x)$ a une singularit\'e isol\'ee \`a l'origine.  Il sera commode d'identifier les modules $M_0$ et $M_1$ \`a $(R\hat\otimes A)^d$.  Rappelons que $A$ d\'esigne un anneau local noetherien complet de corps r\'esiduel $k$. Par MF on entend dans cette section une factorisation matricielle de potentiel $\omega \in \mathcal M_A$, l'id\'eal maximal de $R\hat\otimes A$.  Une d\'eformation de $\overline\M$ \`a $A$ est une classe de MF,  $\M = (\varphi,\psi)$ au-dessus de $A$ telles que par r\'eduction \`a $k$, $\M$ donne $\overline \M$. 
 
 Appelons automorphisme infinit\'esimal  de la d\'eformation $\M$, un automorphisme  de $\M$  qui par r\'eduction  \`a $k$  induit $ 1_{\overline\M}$. Le groupe des $A$-automorphismes de $R\hat\otimes A$ (changements de variables)  agit   sur les MF d\'efinies sur $A$.  Un tel automorphisme est infinit\'esimal si sa r\'eduction \`a $R$ est l'identit\'e.  Une transformation de jauge (d\'efinie sur $A$) est un \'el\'ement du groupe produit semi-direct  du groupe des automorphismes infinit\'esimaux  par le groupe des changements de variables infinit\'esimaux. 
Alors  $\M_1$ et $\M_2$  sont dans  la  m\^eme classe de d\'eformations  de $\overline\M$ \`a $A$,  s'il existe  une transformation de jauge $f:\M_1\stackrel \sim \to \M_2$. Soit $\Def_{\overline\M} (A)$ l'ensemble des d\'eformations de $\overline\M$ \`a $A$, de sorte que $\Def_\M$ est le foncteur des d\'eformations (infinit\'esimales) de $\overline\M$.   Soit l'anneau de Tjurina $\mathcal O = R/(w,\partial w)$. On d\'esigne  enfin par $\mathcal I$ l'id\'eal form\'e des classes de fonctions $Q$-exactes, c'est \`a dire les $h\in R$  tels que $h.1$ est un cobord.  La th\'eorie des d\'eformations de $\overline\M$ se r\'esume en: 
  
   Le foncteur $\Def_{\overline\M}$ admet une enveloppe verselle. L'espace tangent $t_{\overline\M} = \Def_{\overline\M}(k[\epsilon]) \,\, (\epsilon^2 = 0)$ s'ins\`ere dans une suite exacte
 \begin{equation}\label{tan} 0 \to \Ext^1(\overline\M,\overline\M) \to t_{\overline\M} \to \mathcal I \to 0\end{equation}
 Les obstructions au rel\`evement infinit\'esimal de $\overline\M$  appartiennent au conoyau   de l'application $R/ \mathcal I \to \Ext^0(\overline\M,\overline\M)$.  Si  on  d\'eforme   $\overline\M$ en gardant rigide le potentiel $w$,  l'espace tangent  se r\'eduit au  sous-espace $\Ext^1(\overline\M,\overline\M)$, les obstructions \'etant alors localis\'ees   dans $\Ext^0(\overline\M,\overline\M)$. 

 Si $\overline\M$  est munie d'une structure orthogonale  (resp. symplectique, ....) $\overline b$, on s'int\'eresse   aux d\'eformations qui pr\'eservent cette structure. Quitte \`a remplacer une d\'eformation $\M$  par une d\'eformation \'equivalente, on peut supposer que la  structure bilin\'eaire est constante, donn\'ee  par $\overline b\otimes 1$ (resp. $\overline q\otimes 1$). Les isomorphismes (transformations de jauge)  sont  alors assujettis \`a pr\'eserver cette structure. On d\'esignera par un indice sup\'erieur $+$ (resp $-$) le sous-espace des \'el\'ements autoadjoints (resp. anti-autoadjoints). 
 \begin{prop}Le foncteur $\Def_{\overline\M,b}$ (resp. $\Def_{\overline\M,q}$) admet une enveloppe verselle. L'espace tangent $t_{\overline\M,b}  $ (resp. $t_{\overline\M,q}$) est la partie anti-autoadjointe (resp. autoadjointe) de $t_{\overline\M}$, il s'ins\`ere dans une suite exacte
 $0 \to \Ext^1(\overline\M,\overline\M)^{-} \to t_{\overline\M,b} \to \mathcal I \to 0 $ (resp.  $0 \to \Ext^1(\overline\M,\overline\M)^{+} \to t_{\overline\M,q} \to \mathcal I \to 0 $).
 
 Les obstructions au rel\`evement infinit\'esimal de $\overline\M$  appartiennent au conoyau   de l'application $R/ \mathcal I \to \Ext^0(\overline\M,\overline\M)^+$ (resp. $R/ \mathcal I \to \Ext^0(\overline\M,\overline\M)^-)$. \end{prop}

Si  on d\'eforme   $(\overline\M,b)$ (resp. $(\overline\M,q)$) en gardant    $w$ rigide,  l'espace tangent  se r\'eduit au  sous-espace $\Ext^1(\overline\M,\overline\M)^-$ (resp. $\Ext^1(\overline\M,\overline\M)^+$), les obstructions \'etant alors localis\'ees   dans $\Ext^0(\overline\M,\overline\M)^+$ (resp. $\Ext^0(\overline\M,\overline\M)^-)$.

\emph{\underline{ Exemple 3}:} L'exemple cl\'e pour \'enoncer le th\'eor\`eme de p\'eriodicit\'e de Faltings-Kn\"orrer est celui du potentiel $w = xy\in k[[x,y]]$ (point double), et   $\overline\M =  \M_{p,q}$ (exemple 1). On  se limitera \`a  $p=q = r$, avec pour notation $\overline\M_r$.  On a $\overline\M_r  = \overline\M_r^T = \overline\M_r^*$.   On peut prouver qu'\`a    une transformation de jauge pr\`es  on peut supposer que la structure bilin\'eaire (resp. bilin\'eaire tordue) est constante, c'est \`a dire de la forme
 $ b_0 = \left(\begin{array}{cc}  q_0 & 0\\0 &  q_1\end{array}\right),\, \, b_1 = -b_0$, 
avec $q_0,q_1$ d\'efinissant des formes quadratiques non d\'eg\'en\'er\'ees sur $k^r$. M\^eme chose  dans le cas tordu.  

 La d\'eformation verselle de $\overline\M_r$ (resp. $(\overline\M_r,b), \,(\overline\M_r,q)$) est   $\M_{r,{\rm ver}} = \theta(P,Q)$   avec $R_{\rm ver} = W(k)[[p_{ij},q_{ij},t]]/(\star), \,\, (1\leq i,j\leq r)$, ($\star$) d\'esignant l'id\'eal  d\'efini par les \'equations matricielles $(p_{ij}).(q_{ij}) = (q_{ij}).(p_{ij}) = t.1_r$, avec $P = (p_{ij}), \, Q = (q_{ij})$. Dans le cas orthogonal (resp. ...) on impose   $Q = P^{\adj}$ (resp. $Q =  -P^{\adj}$ dans le cas tordu).

 Soit   $(P,Q)$  une MF de rang $r$ sur $A$, de potentiel $\pi \in \mathcal M_A, \,\, (PQ = QP = \pi.1)$.  Formons la MF $(\varphi,\psi)$ de potentiel $xy - \pi$, d\'efinie sur $A[[x,y]]$, par $(\varphi,\psi) = \theta(P,Q) = (Q,-P)\otimes (x,y)$. A une \'equivalence pr\`es  $\varphi = \left(\begin{array}{cc} x  & P\\Q & y \end{array}\right), \,\, \psi = \left(\begin{array}{cc}  y & -P\\-Q & x\end{array}\right)$. Noter que le foncteur $\theta: \MF_A(\pi) \to \MF_{A[[x,y]]}(xy-\pi)$ s'\'etend   aux cat\'egories  triangul\'ees. Si $b $ (resp. $q$) est une forme quadratique, ou symplectique sur $(P,Q)$ (resp. forme quadratique, ou symplectique  tordue), alors $b\otimes 1$ (resp. $q\otimes 1$) est une forme symplectique, ou quadratique  tordue (resp. forme quadratique ou symplectique) sur $\theta(P,Q)$.  On a  (comparer avec \cite{fal})
\begin{thm} i) Toute d\'eformation de $\overline\M_r$ \`a $A$  est isomorphe \`a une d\'eformation de la forme  $\theta(P,Q)$.  Toute d\'eformation de $(M_r,b)$ (resp. $(M_r,q)$) \`a $A$ est isomorphe \`a  $\theta(P,Q)$, pour une  certaine MF orthogonale (symplectique, .....) $(P,Q)$.\\
ii) Le foncteur $\theta: {\rm MF}_A(\pi) \to {\rm MF}_{A[[x,y]]}(xy-\pi)$ est une \'equivalence. \end{thm} 
 Si $A = k[[x_1,\cdots,x_n]]$, et  $\pi\in A\,\,(\pi(0) = 0)$,   $\theta$  \'etablit une \'equivalence ${\rm MF}(\pi) \stackrel\sim \longrightarrow {\rm MF}(xy-\pi)$, resp. ${\rm MF}^{\omega,\pm}(\pi) \stackrel\sim \longrightarrow {\rm MF}^\pm(xy-\pi)$, \,\, ${\rm MF}\pm(\pi) \stackrel\sim \longrightarrow {\rm MF}^{\omega,-(\pm)}(xy-\pi)$.  Le carr\'e  $\theta\theta $  \'echange les  formes quadratiques et symplectiques, et conduit  au r\'esultat de p\'eriodicit\'e d'ordre quatre (comparer avec \cite{orii}):
\begin{thm} Le foncteur $\theta\theta$ d\'efinit une \'equivalence de cat\'egories triangul\'ees
\begin{equation}\label{perio}{\rm MF}^\pm(\pi) \stackrel\sim\longrightarrow {\rm MF}^{-(\pm)} (xy+uv-\pi)\end{equation}\end{thm}
 Lorsque $w$ est une forme quadratique non d\'g\'en\'er\'ee sur un corps $k$, non n\'ec\'essairement alg\'ebriquement clos, de caract\'eristique $\ne 2$, le th\'eor\`eme 6 redonne le fait que la cat\'egorie ${\rm MF}(w)$ est \'equivalente \`a la cat\'egorie des modules de Clifford sur l'alg\`e bre de Clifford de $w$ \cite{ber}, avec, ou sans structure quadratique ou symplectique.


\end{document}